\documentclass[11pt]{article}
\usepackage{pxfonts}
\usepackage{yfonts}
\usepackage{dsfont}
\usepackage{marvosym}
\usepackage{graphicx}
\usepackage{relsize}

\def\bel{\begin{equation}\label}
\def\eeq{\end{equation}}
\def\ds{\displaystyle}
\def\endproof{\hphantom{MM}
\hfill\llap{$\square$}\goodbreak}

\def\R{\mathbb R}

\def\C{\mathfrak{C}}

\def\L{{\bf L}}

\def\I{{\bf I}}
\def\II{{\bf II}}
\def\III{{\bf III}}
\def\IV{{\bf IV}}
\def\M{{\bf M}}
\def\G{{\bf G}}

\def\alpha{\alphaup}
\def\beta{\betaup}
\def\gamma{\gammaup}
\def\delta{\deltaup}
\def\xi{{\xiup}}
\def\eta{{\etaup}}
\def\tau{{\tauup}}
\def\rho{{\rhoup}}
\def\phi{{\phiup}}
\def\psi{{\psiup}}
\def\lambda{{\lambdaup}}
\def\omega{\omegaup}
\def\varphi{{\varphiup}}
\def\gamma{{\gammaup}}

\begin{document}
 \[\begin{array}{cc}\hbox{\LARGE{\bf Fractional integration with singularity}}
 \\\\
  \hbox{\LARGE{\bf on light-cone, I: in $\R^2$}}
 \end{array}\]
 
  \[\hbox{Zipeng Wang}\]
  \[\begin{array}{ccc}
 \hbox{ Department of Mathematics, Westlake university}\\
 \hbox{Cloud town, Hangzhou of China}
 \end{array}\]

 \section{Introduction}
 \setcounter{equation}{0}
 Let $0<\alpha<1$. We consider the fractional integral operator $\I_\alpha$ defined by
 \bel{I_alpha}
 \Big(\I_{\alpha}f\Big)(x,y)~=~\iint_{\R^2} f(x-\mu,y-\nu)\left({1\over |\mu+\nu|}\right)^{1-\alpha}\left({1\over |\mu-\nu|}\right)^{1-\alpha}d\mu d\nu
 \eeq
whose kernel has singularity on the light-cone in $\R^2$.

Throughout the paragraph, we consider $f\ge0$ to be a Schwartz function and denote $\C$ as a generic constant with subindices indicating its dependence.  

{\bf Theorem One~~}   {\it   $\I_\alpha$ defined in (\ref{I_alpha}) extends to a bounded operator 
\bel{Result}
\left\| \I_\alpha f\right\|_{\L^q(\R^2)}~\leq~\C_{p~q}~\left\| f\right\|_{\L^p(\R^2)}, \qquad 1<p<q<\infty
\eeq
if and only if
\bel{Formula}
\alpha~=~{1\over p}-{1\over q}.
\eeq
}

\section{Some preliminary estimates}
\setcounter{equation}{0}
By changing variables $ u={\mu+\nu\over 2}, v={\mu-\nu\over 2}$ in (\ref{I_alpha}), we have
\bel{I_alpha'}
 \Big(\I_{\alpha}f\Big)(x,y)~=~2^{2\alpha-1}\iint_{\R^2} f(x-u-v,y-u+v)\left({1\over |u|}\right)^{1-\alpha}\left({1\over |v|}\right)^{1-\alpha}dudv.
 \eeq

We consider the strong maximal function operator $\M$  defined by 
 \bel{M}
 \Big(\M f\Big)(x,y)~=~\sup_{\gamma>0,~\delta>0}~ {1\over \gamma\delta}\iint_{|u|\leq\delta, ~|v|\leq \gamma} f(x-u-v,y-u+v)dudv
 \eeq
 where the supremum is taking over all rectangles whose  sides are parallel to the lines of degree ${\pi\over 4}$ and ${3\pi\over 4}$. 
 
In particular,  
 \bel{M_1}
  \Big(\M_{\pi\over 4} f\Big)(x,y)~=~\sup_{\delta>0}~ {1\over \delta}\int_{|u|\leq\delta} f(x-u,y-u)du
 \eeq
 is a standard one-parameter maximal function operator defined on the 
  line of degree ${\pi\over 4}$ passing through $(x,y)\in\R^2$.   
  
Vice versa,
  \bel{M_2}
  \Big(\M_{3\pi\over 4} f\Big)(x,y)~=~\sup_{\gamma>0}~ {1\over \gamma}\int_{|v|\leq\gamma} f(x-v,y+v)dv
 \eeq 
 is a standard one-parameter maximal function operator defined on the 
  line of degree ${3\pi\over 4}$ passing through $(x,y)\in\R^2$.

Moreover, by definition given in (\ref{M}),  it is well known
\bel{M<M_1M_2}
\Big(\M f\Big)(x,y)~\leq~\Big(\M_{\pi\over 4} \M_{3\pi\over 4} f\Big)(x,y).
\eeq
Hence that  $\M$ is bounded on $\L^p(\R^2)$ for $1<p<\infty$ whereas
\bel{M_1M_2L^pbound}
\left\| \M_{\pi\over 4} f\right\|_{\L^p(\R^2)}~\leq~\C_p~\left\| f\right\|_{\L^p(\R^2)},\qquad
\left\| \M_{3\pi\over 4} f\right\|_{\L^p(\R^2)}~\leq~\C_p~\left\| f\right\|_{\L^p(\R^2)},\qquad 1<p<\infty.
\eeq
Observe that the norm inequalities in (\ref{M_1M_2L^pbound}) are rotation free. We can choose a new frame whose first coordinate  meets the line of degree ${\pi\over 4}$  or ${3\pi\over 4}$  passing through the origin. 
The $\L^p$-boundedness then follows by applying iteration estimates on every coordinate subspace, as suggested 
  in Chapter II, {\bf 4.1} of the book by Stein \cite{Stein}.

On the other hand, define
 \bel{G}
 \Big(\G f\Big)(x,y)~=~ \left\|\Big(\M_{\pi\over 4} f\Big)(x,y)\right\|_{\L^p(\R)} \left\|\Big(\M_{3\pi\over 4} f\Big)(x,y)\right\|_{\L^p(\R)}
 \eeq 
where
\bel{M_1L^pnorm}
 \left\|\Big(\M_{\pi\over 4} f\Big)(x,y)\right\|_{\L^p(\R)}~=~\left\{\int_{\R} \Big(\M_{\pi\over 4} f\Big)^p(x-v,y+v)dv\right\}^{1\over p}
\eeq
and
\bel{M_2L^pnorm}
  \left\|\Big(\M_{3\pi\over 4} f\Big)(x,y)\right\|_{\L^p(\R)}~=~\left\{\int_{\R} \Big(\M_{3\pi\over 4} f\Big)^p(x-u,y-u)du\right\}^{1\over p}.
\eeq

Let $z={x+y\over 2}, w={x-y\over 2}$. We thus have
\bel{G norm}
\begin{array}{lr}\ds
\iint_{\R^2} \Big(\G f\Big)^p(x,y)dxdy
\\\\ \ds
~=~\iint_{\R^2}\left\{\int_{\R} \Big(\M_{\pi\over 4} f\Big)^p(x-v,y+v)dv\right\}\left\{\int_{\R} \Big(\M_{3\pi\over 4} f\Big)^p(x-u,y-u)du\right\} dxdy
\\\\ \ds
~=~2\iint_{\R^2}\left\{\int_{\R} \Big(\M_{\pi\over 4} f\Big)^p(z+w-v,z-w+v)dv\right\}\left\{\int_{\R} \Big(\M_{3\pi\over 4} f\Big)^p(z+w-u,z-w-u)du\right\} dzdw
\\\\ \ds
~=~2
\left\{\int_{\R}\left\{\int_{\R}\Big(\M_{\pi\over 4} f\Big)^p(z+w-v,z-w+v)dz\right\}dv\right\}
\\\\ \ds~~~~~~~~
\left\{\int\left\{\int_{\R} \Big(\M_{3\pi\over 4} f\Big)^p(z-u+w,z-u-w)dw\right\}du\right\}
\\\\ \ds
~\leq~\C_p~ \left\{\iint_{\R^2}  f^p(z+w-v,z-w+v)dzdv\right\}
\left\{\iint_{\R^2} f^p(z-u+w,z-u-w)dwdu\right\}
\\\\ \ds
~=~\C_p~\left\{\iint_{\R^2} f^p(x,y)dxdy\right\}^2.
\end{array}
\eeq

\section{Proof of Theorem One} 
\setcounter{equation}{0}
Recall from (\ref{I_alpha'}). Given $(x,y)\in\R\times\R$, we consider 
\bel{I_alpha'}
\begin{array}{lr}\ds
\iint_{\R^2} f(x-u-v,y-u+v)\left({1\over |u|}\right)^{1-\alpha}\left({1\over |v|}\right)^{1-\alpha}dudv
 \\\\ \ds
~=~\iint_{\I_{\tau\lambda}\cup\II_{\tau\lambda}\cup\III_{\tau\lambda}\cup\IV_{\tau\lambda}} f(x-u-v,y-u+v)\left({1\over |u|}\right)^{1-\alpha}\left({1\over |v|}\right)^{1-\alpha}dudv
 \end{array}
 \eeq
where
\bel{tau,lambda}
0~<~\tau=\tau(x,y)~<~\infty,\qquad 0~<~\lambda=\lambda(x,y)~<~\infty
\eeq
and
\bel{Regions}
\begin{array}{lr}\ds
\I_{\tau\lambda}~=~\I_{\tau\lambda}(x,y)~=~\left\{|u|\leq \tau, |v|\leq \lambda\right\},\qquad \II_{\tau\lambda}~=~\II_{\tau\lambda}(x,y)~=~\left\{~|u|\leq \tau, |v|> \lambda~\right\},
\\\\ \ds
\III_{\tau\lambda}~=~\III_{\tau\lambda}(x,y)~=~\left\{~|u|>\tau, |v|\leq \lambda~\right\},\qquad
\IV_{\tau\lambda}~=~\IV_{\tau\lambda}(x,y)~=~\left\{~|u|> \tau, |v|> \lambda~\right\}.
\end{array}
\eeq

{\bf 1.} Consider
\bel{norm1}
\begin{array}{lr}\ds
\iint_{\I_{\tau\lambda}} f(x-u-v,y-u+v)\left({1\over |u|}\right)^{1-\alpha}\left({1\over |v|}\right)^{1-\alpha} dudv.
\end{array}
\eeq
We have
\bel{norm1 est}
\begin{array}{lr}\ds
\iint_{\I_{\tau\lambda}} f(x-u-v,y-u+v)\left({1\over |u|}\right)^{1-\alpha}\left({1\over |v|}\right)^{1-\alpha} dudv
\\\\ \ds
~=~\sum_{j=0}^\infty\sum_{k=0}^\infty\iint_{2^{-j-1}\tau<|u|\leq 2^{-j}\tau,~2^{-k-1}\lambda<|v|\leq 2^{-k}\lambda} f(x-u-v,y-u+v)\left({1\over |u|}\right)^{1-\alpha}\left({1\over |v|}\right)^{1-\alpha} dudv
\\\\ \ds
~\leq~\sum_{j=0}^\infty\sum_{k=0}^\infty \left(2^{-j-1}\tau\right)^{\alpha-1}\left(2^{-k-1}\lambda\right)^{\alpha-1}\iint_{2^{-j-1}\tau<|u|\leq 2^{-j}\tau,~2^{-k-1}\lambda<|v|\leq 2^{-k}\lambda} f(x-u-v,y-u+v) dudv
\\\\ \ds
~\leq~\sum_{j=0}^\infty\sum_{k=0}^\infty \left(2^{-j-1}\tau\right)^{\alpha}\left(2^{-k-1}\lambda\right)^{\alpha}
\left\{{4\over 2^{-j}\tau2^{-k}\lambda}\iint_{|u|\leq 2^{-j}\tau, ~|v|\leq 2^{-k}\lambda} f(x-u-v,y-u+v) dudv\right\}
\\\\ \ds
~\leq~4\sum_{j=0}^\infty\sum_{k=0}^\infty \left(2^{-j-1}\right)^{\alpha}\left(2^{-k-1}\right)^{\alpha}~\tau^\alpha\lambda^\alpha \Big(\M f\Big)(x,y)
\\\\ \ds
~\leq~\C_\alpha~\tau^\alpha\lambda^\alpha \Big(\M f\Big)(x,y).
\end{array}
\eeq

{\bf 2.} Consider
\bel{norm2}
\begin{array}{lr}\ds
\iint_{\IV_{\tau\lambda}} f(x-u-v,y-u+v)\left({1\over |u|}\right)^{1-\alpha}\left({1\over |v|}\right)^{1-\alpha} dudv.
\end{array}
\eeq
By applying H\"{o}lder inequality, we have
\bel{norm2 est}
\begin{array}{lr}\ds
 \iint_{\IV_{\tau\lambda}} f(x-u-v,y-u+v)\left({1\over |u|}\right)^{1-\alpha}\left({1\over |v|}\right)^{1-\alpha} dudv
 \\\\ \ds
 ~\leq~ \left\|f\right\|_{\L^p\left(\R^2\right)}
 \left\{\iint_{\IV_{\rho\lambda}} \left({1\over |u|}\right)^{(1-\alpha)\left({p\over p-1}\right)}
 \left({1\over |v|}\right)^{(1-\alpha)\left({p\over p-1}\right)}dudv\right\}^{p-1\over p}.
 \end{array}
\eeq
Note that  $(1-\alpha)\left({p\over p-1}\right)<1$ whereas
\bel{1 p/q>0}
{1 \over p-1}-{\alpha p\over p-1}~=~\left({1\over p}-\alpha\right)\left({p\over p-1}\right)~=~{1\over q}\left({p\over p-1}\right)~>~0
\eeq
by (\ref{Formula}).

We have
\bel{norm2 Est1}
\begin{array}{lr}\ds
\iint_{\IV_{\tau\lambda}}\left({1\over |u|}\right)^{(1-\alpha)\left({p\over p-1}\right)}\left({1\over |v|}\right)^{(1-\alpha)\left({p\over p-1}\right)} dudv
\\\\ \ds
~=~\left\{\int_{\tau\leq|u|} \left({1\over |u|}\right)^{(1-\alpha)\left({p\over p-1}\right)}du\right\}\left\{ \int_{\lambda\leq|v|}\left({1\over |v|}\right)^{(1-\alpha)\left({p\over p-1}\right)}  dv\right\}
\\\\ \ds
~\leq~\C_\alpha~\tau^{\left(\alpha-{1\over p}\right)\left({p\over p-1}\right)}\lambda^{\left(\alpha-{1\over p}\right)\left({p\over p-1}\right)}.
\end{array}
\eeq
By putting together (\ref{norm2 est}) and (\ref{norm2 Est1}), we find
\bel{norm2 L^p}
\begin{array}{lr}\ds
 \iint_{\IV_{\tau\lambda}} f(x-u-v,y-u+v)\left({1\over |u|}\right)^{1-\alpha}\left({1\over |v|}\right)^{1-\alpha} dudv
\\\\ \ds
~\leq~\C_\alpha~ \tau^{\alpha-{1\over p}}\lambda^{\alpha-{1\over p}}\left\|f\right\|_{\L^p(\R^2)}.
\end{array}
\eeq

{\bf 3.} Consider
\bel{norm3}
\iint_{\II_{\tau\lambda}} f(x-u-v,y-u+v)\left({1\over |u|}\right)^{1-\alpha}\left({1\over |v|}\right)^{1-\alpha} dudv.
\eeq
We have
\bel{norm3 est1}
\begin{array}{lr}\ds
\iint_{\II_{\tau\lambda}} f(x-u-v,y-u+v)\left({1\over |u|}\right)^{1-\alpha}\left({1\over |v|}\right)^{1-\alpha} dudv
\\\\ \ds
~=~\int_{|v|>\lambda} \left({1\over |v|}\right)^{1-\alpha} \left\{\sum_{j=0}^\infty \int_{2^{-j-1}\tau<|u|\leq 2^{-j}\tau}
f(x-u-v,y-u+v)\left({1\over |u|}\right)^{1-\alpha} du\right\} dv
\\\\ \ds
~\leq~\int_{|v|>\lambda} \left({1\over |v|}\right)^{1-\alpha} \left\{\sum_{j=0}^\infty \left(2^{-j-1}\tau\right)^{\alpha-1}\int_{2^{-j-1}\tau<|u|\leq 2^{-j}\tau} f(x-u-v,y-u+v) du
\right\} dv
\\\\ \ds
~\leq~\int_{|v|>\lambda} \left({1\over |v|}\right)^{1-\alpha} \left\{\sum_{j=0}^\infty \left(2^{-j-1}\tau\right)^{\alpha}
{4\over 2^{-j}\tau}\int_{|u|\leq 2^{-j}\tau} f(x-u-v,y-u+v) du
\right\} dv
\\\\ \ds
~\leq~\int_{|v|>\lambda} \left({1\over |v|}\right)^{1-\alpha} \left\{4\sum_{j=0}^\infty \left(2^{-j-1}\tau\right)^{\alpha}\right\}
\Big(\M_{\pi\over 4} f\Big)(x-v,y+v) dv
\\\\ \ds
~\leq~\C_\alpha~\tau^{\alpha} \int_{|v|>\lambda} \Big(\M_{\pi\over 4} f\Big)(x-v,y+v) \left({1\over |v|}\right)^{1-\alpha} dv
\end{array}
\eeq
where $\M_{\pi\over 4}$ is defined in (\ref{M_1}).

By applying H\"{o}lder inequality, we have
\bel{norm3 est2}
 \begin{array}{lr}\ds
\int_{|v|>\lambda} \Big(\M_{\pi\over 4} f\Big)(x-v,y+v) \left({1\over |v|}\right)^{1-\alpha} dv
 \\\\ \ds 
 ~\leq~\left\| \Big(\M_{\pi\over 4} f\Big)(x,y)\right\|_{\L^p(\R)}
 \left\{\int_{|v|>\lambda} \left({1\over |v|}\right)^{(1-\alpha)\left({p\over p-1}\right)}dv\right\}^{p-1\over p}
 \\\\ \ds
 ~\leq~ \C_\alpha~ \lambda^{\alpha-{1\over p}}\left\| \Big(\M_{\pi\over 4} f\Big)(x,y)\right\|_{\L^p(\R)} \qquad \hbox{\small{by (\ref{1 p/q>0})}}
  \end{array}
 \eeq
 where $\left\| \Big(\M_{\pi\over 4} f\Big)(x,y)\right\|_{\L^p(\R)}$ is given in (\ref{M_1L^pnorm}).

From (\ref{norm3 est1})-(\ref{norm3 est2}), we have
\bel{norm3 L^p}
\begin{array}{lr}\ds
 \iint_{\II_{\tau\lambda}} f(x-u-v,y-u+v)\left({1\over |u|}\right)^{1-\alpha}\left({1\over |v|}\right)^{1-\alpha} dudv
 \\\\ \ds
 ~\leq~\C_\alpha~\tau^{\alpha} \lambda^{\alpha-{1\over p}}\left\| \Big(\M_{\pi\over 4} f\Big)(x,y)\right\|_{\L^p(\R)}.
\end{array}
\eeq

{\bf 4.} Consider
\bel{norm4}
\iint_{\III_{\tau\lambda}} f(x-u-v,y-u+v)\left({1\over |u|}\right)^{1-\alpha}\left({1\over |v|}\right)^{1-\alpha} dudv.
\eeq
We have
\bel{norm4 est1}
\begin{array}{lr}\ds
\iint_{\III_{\tau\lambda}} f(x-u-v,y-u+v)\left({1\over |u|}\right)^{1-\alpha}\left({1\over |v|}\right)^{1-\alpha} dudv
\\\\ \ds
~=~\int_{|u|>\tau} \left({1\over |u|}\right)^{1-\alpha} \left\{\sum_{k=0}^\infty \int_{2^{-k-1}\lambda<|v|\leq 2^{-k}\lambda}
f(x-u-v,y-u+v)\left({1\over |v|}\right)^{1-\alpha} dv\right\} du
\\\\ \ds
~\leq~\int_{|u|>\tau} \left({1\over |u|}\right)^{1-\alpha} \left\{\sum_{k=0}^\infty \left(2^{-k-1}\lambda\right)^{\alpha-1}\int_{2^{-k-1}\lambda<|v|\leq 2^{-k}\lambda} f(x-u-v,y-u+v) dv
\right\} du
\\\\ \ds
~\leq~\int_{|u|>\tau} \left({1\over |u|}\right)^{1-\alpha} \left\{\sum_{k=0}^\infty \left(2^{-k-1}\lambda\right)^{\alpha}
{4\over 2^{-k}\lambda}\int_{|v|\leq 2^{-k}\lambda} f(x-u-v,y-u+v) dv
\right\} du
\\\\ \ds
~\leq~\int_{|u|>\tau} \left({1\over |u|}\right)^{1-\alpha} \left\{4\sum_{k=0}^\infty \left(2^{-k-1}\lambda\right)^{\alpha}\right\}
\Big(\M_{3\pi\over 4} f\Big)(x-u,y-u) du
\\\\ \ds
~\leq~\C_\alpha~\lambda^{\alpha} \int_{|u|>\tau} \Big(\M_{3\pi\over 4} f\Big)(x-u,y-u) \left({1\over |u|}\right)^{1-\alpha} du
\end{array}
\eeq
where $\M_{3\pi\over 4}$ is defined in (\ref{M_2}).

By applying H\"{o}lder inequality, we have
\bel{norm4 est2}
 \begin{array}{lr}\ds
\int_{|u|>\tau} \Big(\M_{3\pi\over 4} f\Big)(x-u,y-u) \left({1\over |u|}\right)^{1-\alpha} du
 \\\\ \ds 
 ~\leq~\left\| \Big(\M_{3\pi\over 4} f\Big)(x,y)\right\|_{\L^p(\R)}
 \left\{\int_{|u|>\tau} \left({1\over |u|}\right)^{(1-\alpha)\left({p\over p-1}\right)}du\right\}^{p-1\over p}
 \\\\ \ds
 ~\leq~ \C_\alpha~ \tau^{\alpha-{1\over p}}\left\| \Big(\M_{3\pi\over 4} f\Big)(x,y)\right\|_{\L^p(\R)} \qquad \hbox{\small{by (\ref{1 p/q>0})}}
  \end{array}
 \eeq
 where $\left\| \Big(\M_{3\pi\over 4} f\Big)(x,y)\right\|_{\L^p(\R)}$ is given in (\ref{M_2L^pnorm}).

From (\ref{norm4 est1})-(\ref{norm4 est2}), we have
\bel{norm4 L^p}
\begin{array}{lr}\ds
 \iint_{\III_{\tau\lambda}} f(x-u-v,y-u+v)\left({1\over |u|}\right)^{1-\alpha}\left({1\over |v|}\right)^{1-\alpha} dudv
 \\\\ \ds
 ~\leq~\C_\alpha~\tau^{\alpha-{1\over p}} \lambda^{\alpha}\left\| \Big(\M_{3\pi\over 4} f\Big)(x,y)\right\|_{\L^p(\R)}.
\end{array}
\eeq

{\bf 5.} Recall $\G f$ defined in (\ref{G}). Suppose
\bel{Case1}
\Big(\G f\Big)(x,y)~\leq~\Big(\M f\Big)(x,y)\left\|f\right\|_{\L^p\left(\R^2\right)}.
\eeq
We choose $\tau$ and $\lambda$ simultaneously satisfying
\bel{EstCase1}
{\Big(\M f\Big)(x,y)\over \left\| f\right\|_{\L^p\left(\R^2\right)}}~=~\tau^{-{1\over p}}\lambda^{-{1\over p}}
\qquad\hbox{and}\qquad
{ \left\| \Big(\M_{\pi\over 4} f\Big)(x,y)\right\|_{\L^p\left(\R\right)}
\over \left\| \Big(\M_{3\pi\over 4} f\Big)(x,y)\right\|_{\L^p\left(\R\right)} }~=~{ \tau^{-{1\over p}}\over \lambda^{-{1\over p}}}.
\eeq
By solving the equations in (\ref{EstCase1}), we find
\bel{a Case1}
\tau^{-{1\over p}}~=~\left\{ {\Big(\M f\Big)(x,y)\over \left\| f\right\|_{\L^p\left(\R^2\right)}} { \left\| \Big(\M_{\pi\over 4} f\Big)(x,y)\right\|_{\L^p\left(\R\right)}
\over \left\| \Big(\M_{3\pi\over 4} f\Big)(x,y)\right\|_{\L^p\left(\R\right)} }\right\}^{1\over 2}
\eeq
and
\bel{b Case1}
\lambda^{-{1\over p}}~=~\left\{ {\Big(\M f\Big)(x,y)\over \left\| f\right\|_{\L^p\left(\R^2\right)}} { \left\| \Big(\M_{3\pi\over 4} f\Big)(x,y)\right\|_{\L^p\left(\R\right)}
\over \left\| \Big(\M_{\pi\over 4} f\Big)(x,y)\right\|_{\L^p\left(\R\right)} }\right\}^{1\over 2}.
\eeq
On the other hand, suppose
\bel{Case2}
\Big(\G f\Big)(x,y)~>~\Big(\M f\Big)(x,y)\left\|f\right\|_{\L^p\left(\R^2\right)}.
\eeq
We choose $\tau$ and $\lambda$ simultaneously  satisfying
\bel{EstCase2}
{\Big(\G f\Big)(x,y)\over \left\| f\right\|^2_{\L^p\left(\R^2\right)}}~=~\tau^{-{1\over p}}\lambda^{-{1\over p}}
\qquad\hbox{and}\qquad
{ \left\| \Big(\M_{\pi\over 4} f\Big)(x,y)\right\|_{\L^p\left(\R\right)}
\over \left\| \Big(\M_{3\pi\over 4} f\Big)(x,y)\right\|_{\L^p\left(\R\right)} }~=~{ \tau^{-{1\over p}}\over \lambda^{-{1\over p}}}.
\eeq
By solving the equations in (\ref{EstCase2}), we find
\bel{a Case2}
\tau^{-{1\over p}}~=~\left\{ {\Big(\G f\Big)(x,y)\over \left\| f\right\|^2_{\L^p\left(\R^2\right)}} { \left\| \Big(\M_{\pi\over 4} f\Big)(x,y)\right\|_{\L^p\left(\R\right)}
\over \left\| \Big(\M_{3\pi\over 4} f\Big)(x,y)\right\|_{\L^p\left(\R\right)} }\right\}^{1\over 2}
\eeq
and
\bel{b Case2}
\lambda^{-{1\over p}}~=~\left\{ {\Big(\G f\Big)(x,y)\over \left\| f\right\|^2_{\L^p\left(\R^2\right)}} { \left\| \Big(\M_{3\pi\over 4} f\Big)(x,y)\right\|_{\L^p\left(\R\right)}
\over \left\| \Big(\M_{\pi\over 4} f\Big)(x,y)\right\|_{\L^p\left(\R\right)} }\right\}^{1\over 2}.
\eeq

{\bf 6.} Suppose $\Big(\G f\Big)(x,y)\leq\Big(\M f\Big)(x,y)\left\|f\right\|_{\L^p\left(\R^2\right)}$ as (\ref{Case1}).

By inserting (\ref{a Case1})-(\ref{b Case1}) into (\ref{norm1 est}), we have
\bel{norm1 Result Case1}
\begin{array}{lr}\ds
\iint_{\I_{\tau\lambda}} f(x-u-v,y-u+v)\left({1\over |u|}\right)^{1-\alpha}\left({1\over |v|}\right)^{1-\alpha} dudv
\\\\ \ds
~\leq~\C_\alpha~\tau^\alpha\lambda^\alpha \Big(\M f\Big)(x,y)
\\\\ \ds
~\leq~\C_\alpha~ \Big(\M f\Big)(x,y) \left\{ {\Big(\M f\Big)(x,y)\over \left\| f\right\|_{\L^p(\R^2)}}\right\}^{{p\over q}-1} 
~=~\C_\alpha~\Big(\M f\Big)^{p\over q}(x,y)\left\|f\right\|_{\L^p\left(\R^2\right)}^{1-{p\over q}}.
\end{array}
\eeq

By inserting (\ref{a Case1})-(\ref{b Case1}) into (\ref{norm2 L^p}), we have
\bel{norm2 L^p Result Case1}
\begin{array}{lr}\ds
\iint_{\IV_{\tau\lambda}} f(x-u-v,y-u+v)\left({1\over |u|}\right)^{1-\alpha}\left({1\over |v|}\right)^{1-\alpha} dudv
\\\\ \ds
~\leq~\C_\alpha ~\tau^{\alpha-{1\over p}}\lambda^{\alpha-{1\over p}}\left\|f\right\|_{\L^p\left(\R^2\right)} 
\\\\ \ds
~\leq~ \C_\alpha~\left\{ {\Big(\M f\Big)(x,y)\over \left\| f\right\|_{\L^p(\R^2)}}\right\}^{{p\over q}} 
\left\|f\right\|_{\L^p\left(\R^2\right)}
~=~\C_\alpha~\Big(\M f\Big)^{p\over q}(x,y)\left\|f\right\|_{\L^p\left(\R^2\right)}^{1-{p\over q}}.
\end{array}
\eeq

By inserting (\ref{a Case1})-(\ref{b Case1}) into (\ref{norm3 L^p}), we have
\bel{norm3 L^p Result Case1}
\begin{array}{lr}\ds
 \iint_{\II_{\tau\lambda}} f(x-u-v,y-u+v)\left({1\over |u|}\right)^{1-\alpha}\left({1\over |v|}\right)^{1-\alpha} dudv
 \\\\ \ds
 ~\leq~\C_\alpha~\tau^{\alpha} \lambda^{\alpha-{1\over p}}~\left\|\Big(\M_{\pi\over 4} f\Big)(x,y)\right\|_{\L^p\left(\R\right)}

 ~=~\C_\alpha~\tau^{\alpha-{1\over p}} \lambda^{\alpha-{1\over p}}~\tau^{1\over p}~\left\|\Big(\M_{\pi\over 4} f\Big)(x,y)\right\|_{\L^p\left(\R\right)} 
 \\\\ \ds
 ~=~\C_\alpha~\left\{{\Big(\M f\Big)(x,y)\over \left\|f\right\|_{\L^p(\R^2)}}\right\}^{p\over q} 
 
\left\{ {\left\| f\right\|_{\L^p\left(\R^2\right)}\over\Big(\M f\Big)(x,y) } {  \left\| \Big(\M_{3\pi\over 4} f\Big)(x,y)\right\|_{\L^p\left(\R\right)} \over\left\| \Big(\M_{\pi\over 4} f\Big)(x,y)\right\|_{\L^p\left(\R\right)}
 }\right\}^{1\over 2}
\left\|\Big(\M_{\pi\over 4} f\Big)(x,y)\right\|_{\L^p\left(\R\right)}
  \\\\ \ds
~=~\C_\alpha~\left\{{\Big(\M f\Big)(x,y)\over \left\|f\right\|_{\L^p(\R^2)}}\right\}^{p\over q}
\left\{{\Big(\G f\Big)(x,y)\over \Big( \M f\Big)(x,y)\left\|f\right\|_{\L^p(\R^2)}}\right\}^{1\over 2}
\left\|f\right\|_{\L^p\left(\R^2\right)}
 \\\\ \ds
 ~\leq~\C_\alpha~\Big(\M f\Big)^{p\over q}(x,y)\left\|f\right\|_{\L^p\left(\R^2\right)}^{1-{p\over q}}.
 \end{array}
\eeq
On the other hand, by inserting (\ref{a Case1})-(\ref{b Case1}) into (\ref{norm4 L^p}), we have
\bel{norm4 L^p Result Case1}
\begin{array}{lr}\ds
 \iint_{\III_{\tau\lambda}} f(x-u-v,y-u+v)\left({1\over |u|}\right)^{1-\alpha}\left({1\over |v|}\right)^{1-\alpha} dudv
 \\\\ \ds
 ~\leq~\C_\alpha~\tau^{\alpha-{1\over p}} \lambda^{\alpha}~\left\|\Big(\M_{3\pi\over 4} f\Big)(x,y)\right\|_{\L^p\left(\R\right)}

 ~=~\C_\alpha~\tau^{\alpha-{1\over p}} \lambda^{\alpha-{1\over p}}~\lambda^{1\over p}~\left\|\Big(\M_{3\pi\over 4} f\Big)(x,y)\right\|_{\L^p\left(\R\right)} 
 \\\\ \ds
 ~=~\C_\alpha~\left\{{\Big(\M f\Big)(x,y)\over \left\|f\right\|_{\L^p(\R^2)}}\right\}^{p\over q} 
 
\left\{ {\left\| f\right\|_{\L^p\left(\R^2\right)}\over\Big(\M f\Big)(x,y) } {  \left\| \Big(\M_{\pi\over 4} f\Big)(x,y)\right\|_{\L^p\left(\R\right)} \over\left\| \Big(\M_{3\pi\over 4} f\Big)(x,y)\right\|_{\L^p\left(\R\right)}
 }\right\}^{1\over 2}
\left\|\Big(\M_{3\pi\over 4} f\Big)(x,y)\right\|_{\L^p\left(\R\right)}
  \\\\ \ds
~=~\C_\alpha~\left\{{\Big(\M f\Big)(x,y)\over \left\|f\right\|_{\L^p(\R^2)}}\right\}^{p\over q}
\left\{{\Big(\G f\Big)(x,y)\over \Big( \M f\Big)(x,y)\left\|f\right\|_{\L^p(\R^2)}}\right\}^{1\over 2}
\left\|f\right\|_{\L^p\left(\R^2\right)}
 \\\\ \ds
 ~\leq~\C_\alpha~\Big(\M f\Big)^{p\over q}(x,y)\left\|f\right\|_{\L^p\left(\R^2\right)}^{1-{p\over q}}.
 \end{array}
\eeq

{\bf 7.} Suppose $\Big(\G f\Big)(x,y)>\Big(\M f\Big)(x,y)\left\|f\right\|_{\L^p\left(\R^2\right)}$ as (\ref{Case2}).

By inserting (\ref{a Case2})-(\ref{b Case2}) into (\ref{norm1 est}), we have
\bel{norm1 Result Case2}
\begin{array}{lr}\ds
\iint_{\I_{\tau\lambda}} f(x-u-v,y-u+v)\left({1\over |u|}\right)^{1-\alpha}\left({1\over |v|}\right)^{1-\alpha} dudv
\\\\ \ds
~\leq~\C_\alpha~\tau^\alpha\lambda^\alpha \Big(\M f\Big)(x,y)
\\\\ \ds
~\leq~\C_\alpha~\tau^{\alpha} \lambda^{\beta}\left\{{\Big(\G f\Big)(x,y)\over \left\|f\right\|_{\L^p\left(\R^2\right)}}\right\}
\\\\ \ds
~=~\C_\alpha~\left\{{\Big(\G f\Big)(x,y)\over \left\|f\right\|^2_{\L^p\left(\R^2\right)}}\right\}^{{p\over q}-1} \left\{{\Big(\G f\Big)(x,y)\over \left\|f\right\|_{\L^p\left(\R^2\right)}}\right\}
~=~\C_\alpha~\Big(\G f\Big)^{p\over q}(x,y)\left\|f\right\|_{\L^p\left(\R^2\right)}^{1-{2p\over q}}.
\end{array}
\eeq

By inserting (\ref{a Case2})-(\ref{b Case2}) into (\ref{norm2 L^p}), we have
\bel{norm2 L^p Result Case2}
\begin{array}{lr}\ds
\iint_{\IV_{\tau\lambda}} f(x-u-v,y-u+v)\left({1\over |u|}\right)^{1-\alpha}\left({1\over |v|}\right)^{1-\alpha} dudv
\\\\ \ds
~\leq~\C_\alpha ~\tau^{\alpha-{1\over p}}\lambda^{\alpha-{1\over p}}\left\|f\right\|_{\L^p\left(\R^2\right)} 
\\\\ \ds
~\leq~ \C_\alpha~\left\{ {\Big(\G f\Big)(x,y)\over \left\| f\right\|^2_{\L^p(\R^2)}}\right\}^{{p\over q}} 
\left\|f\right\|_{\L^p\left(\R^2\right)}
~=~\Big(\G f\Big)^{p\over q}(x,y)\left\|f\right\|_{\L^p\left(\R^2\right)}^{1-{2p\over q}}.
\end{array}
\eeq

By inserting (\ref{a Case2})-(\ref{b Case2}) into (\ref{norm3 L^p}), we have
\bel{norm3 L^p Result Case1}
\begin{array}{lr}\ds
 \iint_{\II_{\tau\lambda}} f(x-u-v,y-u+v)\left({1\over |u|}\right)^{1-\alpha}\left({1\over |v|}\right)^{1-\alpha} dudv
 \\\\ \ds
 ~\leq~\C_\alpha~\tau^{\alpha} \lambda^{\alpha-{1\over p}}~\left\|\Big(\M_{\pi\over 4} f\Big)(x,y)\right\|_{\L^p\left(\R\right)}

 ~=~\C_\alpha~\tau^{\alpha-{1\over p}} \lambda^{\alpha-{1\over p}}~\tau^{1\over p}~\left\|\Big(\M_{\pi\over 4} f\Big)(x,y)\right\|_{\L^p\left(\R\right)} 
 \\\\ \ds
 ~=~\C_\alpha~\left\{{\Big(\G f\Big)(x,y)\over \left\|f\right\|^2_{\L^p(\R^2)}}\right\}^{p\over q} 
 
\left\{ {\left\| f\right\|^2_{\L^p\left(\R^2\right)}\over\Big(\G f\Big)(x,y) } {  \left\| \Big(\M_{3\pi\over 4} f\Big)(x,y)\right\|_{\L^p\left(\R\right)} \over\left\| \Big(\M_{\pi\over 4} f\Big)(x,y)\right\|_{\L^p\left(\R\right)}
 }\right\}^{1\over 2}
\left\|\Big(\M_{\pi\over 4} f\Big)(x,y)\right\|_{\L^p\left(\R\right)}
  \\\\ \ds
~=~\C_\alpha~\left\{{\Big(\G f\Big)(x,y)\over \left\|f\right\|^2_{\L^p(\R^2)}}\right\}^{p\over q}
\left\|f\right\|_{\L^p\left(\R^2\right)}
 ~=~\C_\alpha~\Big(\G f\Big)^{p\over q}(x,y)\left\|f\right\|_{\L^p\left(\R^2\right)}^{1-{2p\over q}}.
 \end{array}
\eeq
On the other hand, by inserting
(\ref{a Case2})-(\ref{b Case2}) into (\ref{norm4 L^p}), we have
\bel{norm3 L^p Result Case1}
\begin{array}{lr}\ds
 \iint_{\III_{\tau\lambda}} f(x-u-v,y-u+v)\left({1\over |u|}\right)^{1-\alpha}\left({1\over |v|}\right)^{1-\alpha} dudv
 \\\\ \ds
 ~\leq~\C_\alpha~\tau^{\alpha-{1\over p}} \lambda^{\alpha}~\left\|\Big(\M_{3\pi\over 4} f\Big)(x,y)\right\|_{\L^p\left(\R\right)}

 ~=~\C_\alpha~\tau^{\alpha-{1\over p}} \lambda^{\alpha-{1\over p}}~\lambda^{1\over p}~\left\|\Big(\M_{3\pi\over 4} f\Big)(x,y)\right\|_{\L^p\left(\R\right)} 
 \\\\ \ds
 ~=~\C_\alpha~\left\{{\Big(\G f\Big)(x,y)\over \left\|f\right\|^2_{\L^p(\R^2)}}\right\}^{p\over q} 
 
\left\{ {\left\| f\right\|^2_{\L^p\left(\R^2\right)}\over\Big(\G f\Big)(x,y) } {  \left\| \Big(\M_{\pi\over 4} f\Big)(x,y)\right\|_{\L^p\left(\R\right)} \over\left\| \Big(\M_{3\pi\over 4} f\Big)(x,y)\right\|_{\L^p\left(\R\right)}
 }\right\}^{1\over 2}
\left\|\Big(\M_{3\pi\over 4} f\Big)(x,y)\right\|_{\L^p\left(\R\right)}
  \\\\ \ds
~=~\C_\alpha~\left\{{\Big(\G f\Big)(x,y)\over \left\|f\right\|^2_{\L^p(\R^2)}}\right\}^{p\over q}
\left\|f\right\|_{\L^p\left(\R^2\right)}
~=~\C_\alpha~\Big(\G f\Big)^{p\over q}(x,y)\left\|f\right\|_{\L^p\left(\R^2\right)}^{1-{p\over q}}.
 \end{array}
 \eeq
 \endproof

wangzipeng@westlake.edu.cn

\end{document}